\UseRawInputEncoding
\documentclass[11pt]{article}
\usepackage{amsmath,amsfonts,amsthm,amsopn,amssymb}
\usepackage{cite,marginnote}
\usepackage{color,enumitem,graphicx}
\usepackage[colorlinks=true,urlcolor=blue,
citecolor=red,linkcolor=blue,linktocpage,pdfpagelabels,
bookmarksnumbered,bookmarksopen]{hyperref}
\usepackage[english]{babel}

\usepackage[left=2.9cm,right=2.9cm,top=2.8cm,bottom=2.8cm]{geometry}

\numberwithin{equation}{section}

\makeindex
\makeindex
\newtheorem{theorem}{Theorem}[section]
\newtheorem{definition}[theorem]{Definition}
\newtheorem{lemma}[theorem]{Lemma}
\newtheorem{corollary}[theorem]{Corollary}
\newtheorem{proposition}[theorem]{Proposition}
\newtheorem{remark}[theorem]{Remark}

\newcommand{\s}{\section}

\newcommand{\R}{\mathbb R}

\newcommand{\lab}{\label}
\newcommand{\bt}{\begin{theorem}}
\newcommand{\et}{\end{theorem}}
\newcommand{\bl}{\begin{lemma}}
\newcommand{\el}{\end{lemma}}
\newcommand{\bd}{\begin{definition}}
\newcommand{\ed}{\end{definition}}
\newcommand{\bc}{\begin{corollary}}
\newcommand{\ec}{\end{corollary}}
\newcommand{\bp}{\begin{proof}}
\newcommand{\ep}{\end{proof}}
\newcommand{\bx}{\begin{example}}
\newcommand{\ex}{\end{example}}
\newcommand{\bi}{\begin{exercise}}
\newcommand{\ei}{\end{exercise}}
\newcommand{\bo}{\begin{proposition}}
\newcommand{\eo}{\end{proposition}}
\newcommand{\br}{\begin{remark}}
\newcommand{\er}{\end{remark}}
\newcommand{\beq}{\begin{equation}}
\newcommand{\eeq}{\end{equation}}
\newcommand{\ba}{\begin{align}}
\newcommand{\ea}{\end{align}}
\newcommand{\bn}{\begin{enumerate}}
\newcommand{\en}{\end{enumerate}}
\newcommand{\bg}{\begin{align*}}
\newcommand{\bcs}{\begin{cases}}
\newcommand{\ecs}{\end{cases}}

\newcommand{\NN}{{\mathbb N}}

\newcommand{\bean}{\begin{eqnarray*}}
\newcommand{\eean}{\end{eqnarray*}}

\def\N{\mathbb{N}}

\def\R{\mathbb{R}}

\def\bd{\mathrm{bd}\,}

\begin{document}

\title{Bifurcation from essential spectrum for an elliptic equation with general nonlinearity}
\author{Jianjun Zhang$^a$, Xuexiu Zhong$^b$, Huan-Song Zhou$^c$\thanks{Corresponding author. \newline Email address: J.J.Zhang: zhangjianjun09@tsinghua.org.cn; X.X.Zhong: zhongxuexiu1989@163.com; H.S.Zhou: hszhou@whut.edu.cn. }\\
\vspace{-0.4cm}\\
\small $^a$College of Mathematics and Statistics, Chongqing Jiaotong University,\\
\small Chongqing 400074, PR China\\
\vspace{-0.4cm} \ \\
\small  $^b$South China Research Center for Applied Mathematics and Interdisciplinary Studies,\\
\small South China Normal University, Guangzhou 510631, PR China\\
\vspace{-0.4cm} \ \\
\small $^c$ Center for Mathematical Sciences and Department of Mathematics,\\
\small Wuhan University of Technology, Wuhan 430070, PR China
}





\date{}
\maketitle

\begin{abstract}
In this paper, based on some prior estimates, we show that the essential spectrum $\lambda=0$ is a bifurcation point for an superlinear  elliptic equation with only local conditions, which improves the related results on an open problem proposed by C.A. Stuart in [Lecture Notes in Mathematics, 1017].
\end{abstract}

\indent  {\it Declarations of interest:} none.

\s{Introduction}
\renewcommand{\theequation}{1.\arabic{equation}}
In the present paper, we are concerned with the bifurcation for the following elliptic problem
\beq\lab{eq:B-P}
\begin{cases}
-\Delta u +\lambda u=q(x)g(u)\;&\hbox{in}~\R^N, N\geq 1,\\
(\lambda, u)\in \R\times H^1(\R^N),& u\not\equiv 0,
\end{cases}
\eeq
where $g(s)$ and $q(x)$ satisfy the following assumptions
\begin{itemize}
\item[(G1)] there exists $s_1>0$ such that $g\in C([0,s_1])$,  $g(0)=0$ and $g(s)>0$ for $s\in (0,s_1]$;
\item[(G2)] there exist $\sigma\in (0,\frac{4}{N})$ and $A>0$, such that
$$\lim_{s\rightarrow 0^+}\frac{g(s)}{s^{\sigma+1}}=A.$$
\item[(Q)] $q(x)\in L^\infty(\R^N)$ such that
$$q_0:=\sup_{x\in \R^N}q(x)>\inf_{x\in \R^N}q(x)=\lim_{|x|\rightarrow +\infty}q(x)=:q_\infty>0;$$
\end{itemize}
We say $\lambda =0$ is a bifurcation point for \eqref{eq:B-P}, if there is a sequence $\{ \lambda_n, u_n\}$ of solutions for \eqref{eq:B-P} such that $ \lambda_n \rightarrow 0$ and $\| u_n\|_{H^1} \rightarrow 0$ as $n \rightarrow 0$.


It is well known that the Laplacian operator $-\Delta $ admits no eigenvalues in $H^1(\R^N)$ and has only the essential spectrum $\lambda=0$. Hence the standard Lyapunov-Schmidt reduction is not applicable to problem \eqref{eq:B-P}. It is worth mentioning that $q(x)$ in this paper is not necessary to be radial. Hence, we can not work in the subspace of radial functions. So, even for $\lambda>0$, the operator $(-\Delta+\lambda)^{-1}:H^1(\R^N)\rightarrow H^1(\R^N)$ is not compact. In particular, if $q(x)\not\rightarrow 0$ as $|x|\rightarrow +\infty$, the operator $$\mathbb{A}_\lambda:=(-\Delta+\lambda)^{-1} \circ(q(\cdot)g)$$ is not completely continuous in $H^1(\R^N)$. Thus, it seems that the associated Leray-Schauder degree is not well defined in $H^1(\R^N)$ and the method of topological degree fails.
Hence, the topic about bifurcation from the essential spectrum of nonlinear differential equations, was once a very challenging problem and of general interest in the 1980s, see \cite{BenciFortunato1981,BenciFortunato1981-b,Heinz1986,BHK1983,Jeanjean1999,Stuart1980,Stuart1981,Stuart1982,Stuart1983,ZhuZhou1988,Stuart1989} and references therein.

In paper \cite{BenciFortunato1981},  Benci and Fortunato studied \eqref{eq:B-P} with
\begin{center}
$g(u)=|u|^{p-2}u$,  $-q(x)\geq const \cdot |x|^\alpha$ and  $\alpha>\frac{N(p-2)}{2}$,
\end{center}
and they proved that $\lambda=\lambda_0$ is a bifurcation point for any $\lambda_0\leq 0$, by virtue of certain compact embeddings of a weight Sobolev space.  In \cite{BHK1983}, if $0\leq q(x)\in C(\R^N)\cap L^1(\R^N)$, the authors showed that there are infinitely many global solution ``branches" emanating from the trivial line of solutions at $\lambda^*=0$.
By applying the compactness of radial Sobolev space and Pohozaev identity,    Berestycki and Lions in  \cite{BerestyckiLions1980}  considered problem  \eqref{eq:B-P} in the autonomous case, that is, $q(x) \equiv 1$ in \eqref{eq:B-P}.
Stuart \cite{Stuart1982,Stuart1983} investigated problem \eqref{eq:B-P} with $g(u)=|u|^{\sigma}u$ and $q(x)$ being non-radially symmetric, to obtain the compactness, he imposed a vanishing condition on $q(x)$, that is,
 $$\displaystyle \lim_{|x|\rightarrow +\infty}q(x)=0,$$ which can guarantee that the map $H^1(\R^N)\ni u\mapsto \int_{\R^N}q(x)G(u)dx$ is weakly sequentially continuous. Actually, in \cite{Stuart1982,Stuart1983},  the following $L^2-$normalized problem
\beq\lab{eq:20211020-e1}
\begin{cases}
-\Delta u+\lambda u=q(x)|u|^\sigma u~\hbox{in}~\R^N, \ \sigma\in (0,\frac{4}{N}),\\
\int_{\R^N}|u|^2 dx=c >0 \text{ and }  c\in \mathbb{R},\\
\end{cases}
\eeq
was considered. By further assuming that
\beq\lab{eq:;20220102-e1}
q(x)\geq A(1+|x|)^{-t}, x\in \R^N,\,\,\mbox{for some}\,\, A>0, t\in (0,2-\frac{N\sigma}{2}).
\eeq
 Stuart proved the existence of $L^2-$normalized solution $(\lambda_c,u_c)$ (i.e., $\|u_c\|_{L^2(\R^N)}^{2}=c$) for $c>0$ small. Furthermore, $\lambda_c\rightarrow 0$ in $\R$ and $u_c\rightarrow 0$ in $H^1(\R^N)$ as $c\rightarrow 0^+$, that is, $\lambda=0$ is a bifurcation point for \eqref{eq:B-P}.
We note that the assumption \eqref{eq:;20220102-e1} plays a crucial role in ensuring the positivity of the Lagrange multiplier $\lambda_c>0$ for given $c>0$.
So, Stuart proposed in \cite[Remark 2, Page 577]{Stuart1983}  an open problem as follows
\begin{center}
	{\it ``At present, no results seem to cover cases where  $q(x)\not\rightarrow 0$ as $|x|\rightarrow +\infty$ and $q$ is non-radially symmetric." }
\end{center}
 After that, paper \cite{ZhuZhou1988} gave an answer to the above open problem for  \eqref{eq:B-P}.
Precisely, Zhu and Zhou \cite{ZhuZhou1988} considered a much more general problem
\beq\lab{eq:20210928-e1}
-\Delta u+\lambda u=f(x, u)~\hbox{in}~\R^N, (\lambda, u)\in \R\times H^1(\R^N)
\eeq
and they required  $f(x, u)$ satisfying the following conditions
\begin{itemize}
\item[(f1)]$f(x, u)\in C(\R^N \times \R)$ and for some constant $C>0$ such that
$$|f(x, u)|\leq C(|u|+|u|^{1+\frac{4}{N}}), \forall x\in \R^N, u\in \R;$$
\item[(f2)]$\displaystyle \lim_{u\rightarrow 0}u^{-1}f(x, u)=0$ uniformly in $x\in \R^N$;
\item[(f3)]$f(x, u)\rightarrow \bar{f}(u)$ as $|x|\rightarrow +\infty$ uniformly for bounded $u\in \R$;
\item[(f4)]$\exists A>0, s\in (0,2)$ and $\beta\in (2,2+\frac{2(2-s)}{N})$ such that $$F(x, t)\geq A(1+|x|)^{-s}|t|^\beta, x\in \R^N,t\in \R,$$ or $$\displaystyle \lim_{t\rightarrow 0}\frac{\bar{F}(t)}{|t|^{2+\frac{4}{N}}}=+\infty;$$
\item[(f5)]$F(x,t)\geq \bar{F}(t), \forall x\in \R^N, t\in \R$.
\end{itemize}
Thanks to the concentration compactness principle, in \cite{ZhuZhou1988} they proved that: For any $c>0$ small enough, there is  $\lambda_c>0$ such that problem \eqref{eq:20210928-e1}
possesses a $L^2$-normalized solution $u_c$,  
that is, $(\lambda_c, u_c)$ solves \eqref{eq:20210928-e1} and $I(u_c)=m_c$, where
$$I(u):=\frac{1}{2}\int_{\R^N}|\nabla u|^2 dx-\int_{\R^N}F(x, u)dx, u\in H^1(\R^N)$$
and
$$m_c=\inf\left\{I(u): \|u\|_{L^2(\R^N)}^{2}=c\right\}, \ c>0.$$
Furthermore, they showed that $\lambda_c\rightarrow 0$ in $\R$ and $u_c\rightarrow 0$ in $H^1(\R^N)$ as $c\rightarrow 0^+$, that is, $\lambda=0$ is a bifurcation point for \eqref{eq:20210928-e1}.  Meanwhile, Stuart \cite{Stuart1989} also studied \eqref{eq:20210928-e1} with $f(x, u)=\sum_{i=1}^{m}q_i(x)|u|^{\sigma_i}u$ and $0<\sigma_i<\frac{4}{N}$ and obtained 
the existence of normalized ground state solution $(\lambda_c, u_c)$
with $\lambda_c\rightarrow 0, \|\nabla u_c\|_2^2\rightarrow 0$ as $c\rightarrow 0^+$.
It is known that to get the compactness for a minimizing sequence constrained on $S_c$, the Lagrange multiplier has to be controlled. So, in  \cite{Stuart1989}, Stuart still required the conditions
\begin{center}
$q_i(x)\geq A_i(1+|x|)^{-t_i}$ with $A_i>0, t_i\in [0,2)$ and $0<\sigma_i<\frac{2(2-t_i)}{N}$,
\end{center}
which play a crucial role in verifying $\inf_{u\in S_c}J[u]<0$ (see the conditions (A3) and (A4) in \cite{Stuart1989}). In particular, $q(x)\not\rightarrow 0$ as $|x|\rightarrow+\infty$ is also involved there, see the condition (A4) in \cite{Stuart1989}.
Some progress was also implied in \cite{Shibata2014}, although the author did not give a straightforward estimate on the Lagrange multiplier that $\lambda_c\rightarrow 0$ as $c\rightarrow 0$.

 However, all the papers mentioned above are required some global properties on the nonlinearity $g(s)$.
 We note that Jeanjean in \cite{Jeanjean1999} studied the local conditions on nonlinearities to ensure bifurcation from the continuous spectrum for
\eqref{eq:20210928-e1}, but he still required that $\displaystyle \lim_{|x|\rightarrow +\infty}f(x, s)=0$ uniformly for $s\in [-\delta,\delta]$, see \cite[hypothesis (H2)]{Jeanjean1999}, which is essentially not th-e case of Stuart's open problem.
 To the best of the authors' knowledge, there seemed no further results on Stuart's open problem under more general nonlinear terms. In this paper, our main aim is to answer this open problem by requiring only the properties of $g(s)$ near the origin. 
 Our main theorem is as follows:
\bt\lab{thm:main}
Assume (G1)-(G2) and (Q) hold, then $\lambda=0$ is a bifurcation point for problem \eqref{eq:B-P}.
\et

\br\lab{remark:20220628-r1}
{\bf (i) } The well known Clark's theorem \cite{Clark1973} asserts the existence of a sequence of negative critical values tending to $0$ for even coercive functionals, which has been improved in many works, see, e.g., \cite{Heinz1987,Wang2001,LiuWang2015,CLW2017}.  In particular, Wang \cite{Wang2001} requires only  the local conditions near the origin.
However, these works are focused on the sub-linear case. In this paper, our aim is  to establish the bifurcation results for a super-linear problem with only local conditions near the origin,  the Clark's theorem and its  variants are not applicable.
{\bf (ii)}  We mention that, if $-\Delta$ is replaced by a Schr\"odinger operator $-\Delta +V(x)$ with $V(x)$ being a coercive potential and $f(x,u)$ being a power type nonlinearity in \eqref{eq:20210928-e1}, the $L^2-$normalized solutions and its asymptotical behaviors for \eqref{eq:20210928-e1} has attracted much attentions in recent years, see, e.g.,  \cite{GuoSeiringer2014, GuoZengZhou2016,  LiYangZou2020, YangYang2022}. Very recently, some bifurcation results for a nonhomogeneous fractional equation was obtained in \cite{HeZou2020}.
\er

To prove the above Theorem \ref{thm:main}, we should mention a very recent paper due to Jeanjean-Zhang-Zhong \cite{JeanZhangZhong2021}, in which a new approach is  developed to study the $L^2$-normalized solutions
for the following equation
 \begin{equation}\label{eq:NLS}
-\Delta u+\lambda u=g(u)\;\hbox{in}\;\R^N,
\end{equation}
 under very mild assumptions on $g$. In \cite{JeanZhangZhong2021},
the authors obtained some results about the existence, non-existence and multiplicities of the $L^2$-normalized solutions. In particular, under some suitable assumptions, they established the asymptotic behavior of positive solutions to \eqref{eq:NLS} as $\lambda\rightarrow 0^+$.
In the spirit of the work \cite{JeanZhangZhong2021},  our proof for Theorem \ref{thm:main} is divided into the following steps:
\begin{itemize}
\item[(i)]We firstly focus on the following modified equation of \eqref{eq:B-P}
\beq\lab{eq:20210915-e1}
-\Delta u+\lambda u=q(x)\tilde{g}(u)~\hbox{in}~\R^N.
\eeq
where $\tilde{g}(s)=g(s)$ provided $s$ small enough and satisfies some suitable assumptions, see Lemma \ref{lemma:20210915-l1} for the details.
\item[(ii)]We show that \eqref{eq:20210915-e1} possesses a mountain pass solution $u_\lambda$ for any $\lambda>0$, see our Theorem \ref{thm:mountain-pass} in Section \ref{sec:existence}.
\item[(iii)] We then prove that $u_\lambda\rightarrow 0$ in $H^1(\R^N)$ as $\lambda\rightarrow 0^+$, that is, $\lambda=0$ is a bifurcation point for equation \eqref{eq:20210915-e1}, see our Proposition \ref{prop:nb-p2}.
\item[(iv)] We finally claim that $\|u_\lambda\|_\infty\rightarrow 0$ as $\lambda\rightarrow 0^+$ (see our Lemma \ref{lemma:20210915-l3} or Remark \ref{remark:20220105-br1}), which shows that  $u_\lambda$ is essentially a solution to the original problem \eqref{eq:B-P} for $\lambda$ small enough. Hence, $\lambda=0$ is also a bifurcation point for equation \eqref{eq:B-P}.
\end{itemize}

\br\lab{remark:20211016-r1}
In paper \cite{ZhuZhou1988}, the open problem of \cite{Stuart1983} was answered under the assumptions (f1)-(f5), but (f1) is a global condition on $f(x, u)$.
Condition (A1) in \cite{Stuart1989} is also a global condition on $g(s)$ and $g(s)$ in \cite{Stuart1989} is homogeneous. In this paper, we start with the modified problem \eqref{eq:20210915-e1}, our main results only depends on the local behavior of $g(s)$ near $0$. In our case, $g(s)$ is allowed to grow very fast, such as, $g(s)$ can be mass super-critical or even Sobolev super-critical at infinity. Moreover, $g(s)$ is also allowed to be negative and decreasing for large $s>0$.
\er

\br\lab{remark:20211016-r2}
By the arguments of our paper, it is easy to generalize the result for the following problem
    \beq\lab{eq:20211016-e1}
    -\Delta u+\lambda u=\sum_{i=1}^{m}q_i(x)g_i(u)~\hbox{in}~\R^N,\\
    (\lambda, u)\in \R\times H^1(\R^N).
    \eeq
    Assume that for each $i\in \{1,2,\cdots,m\}$, $g_i$ satisfies (G1)-(G2) and $q_i$ admits (Q), then $\lambda=0$ is also a bifurcation point of \eqref{eq:20211016-e1}.
\er

Based on Theorem \ref{thm:main}, we have also a bifurcation result for problem \eqref{eq:B-P}  with $\lambda u$ being replaced by $\lambda V(x)u$ in \eqref{eq:B-P}, see \eqref{eq:20210917-ze1}.
\bt\lab{cro:20210917-c1}
Let $q(x)\equiv C>0$ or satisfies (Q). If the conditions (G1)-(G2) hold and $V(x)$ satisfies
\[
0<V_0:=\inf_{x\in \R^N}V(x)<\sup_{x\in \R^N}V(x)=\lim_{|x|\rightarrow \infty}V(x)=:V_\infty,
\leqno{(V)}
\]
then $\lambda=0$ is a bifurcation point for the following equation
\beq\lab{eq:20210917-ze1}
-\Delta u+\lambda V(x)u=q(x)g(u), (\lambda, u)\in \R\times H^1(\R^N).
\eeq
\et

\s{Preliminaries}\lab{sec:preliminary}
\renewcommand{\theequation}{2.\arabic{equation}}
\bl\lab{lemma:20210915-l1}
Under the assumptions (G1)-(G2), for any $\alpha,\beta$ satisfying $2<\alpha<2+\sigma<\beta<2+\frac{4}{N}$, there exist $\tilde{g}:\R\rightarrow \R$ and $s_0\in (0,s_1)$ such that
\begin{itemize}
\item[($\tilde{G}1$)] $\tilde{g}\in C(\R)$ is odd and $\tilde{g}(s)>0$ for $s>0$;
\item[($\tilde{G}2$)]$\tilde{g}(s)\equiv g(s)$ for $s\in [0,s_0]$ and $\displaystyle \lim_{s\rightarrow +\infty}\frac{\tilde{g}(s)}{s^{\beta-1}}=:B>0$;
\item[($\tilde{G}3$)]$\displaystyle \alpha \tilde{G}(s)\leq \tilde{g}(s)s\leq \beta\tilde{G}(s)$ for any  $s\in \R$, where $\tilde{G}(s)=\int_0^s \tilde{g}(t)dt$.
\end{itemize}
\el
\bp
Under the assumptions (G1)-(G2), we know that
$$g(s)=(A+o(1))s^{1+\sigma}, G(s)=\frac{1}{2+\sigma}(A+o(1))s^{2+\sigma}\;\hbox{for $s>0$ close to $0$},$$
where we denote by $o(1)$ a quantity which goes to $0$ as $s\rightarrow 0$.
Hence, by $\alpha<2+\sigma<\beta$, there exists $s_0\in (0,s_1)$ such that
\beq\lab{eq:20220102-e2}
\alpha G(s)\leq g(s)s\leq \beta G(s), \forall s\in [-s_0,s_0].
\eeq
We define
\beq\lab{eq:def-tilde-g}
\tilde{g}(s):=
\begin{cases}
g(s),\quad &\hbox{for}~s\in[0,s_0],\\
g(s_0)s_{0}^{-(\beta-1)} s^{\beta-1},\;&\hbox{for}~s\geq s_0,\\
-\tilde{g}(-s),\;&\hbox{for}~s\leq 0.
\end{cases}
\eeq
Clearly, $\tilde{g}\in C(\R)$ is odd and $\tilde{g}(s)>0$ for $s>0$, thus $(\tilde{G}1)$ holds. By the definition \eqref{eq:def-tilde-g}, it is easy to see that $(\tilde{G}2)$ holds with  $$B:=\lim_{s\rightarrow +\infty}\frac{\tilde{g}(s)}{s^{\beta-1}}=g(s_0)s_{0}^{-(\beta-1)}>0.$$

Furthermore, for $s>s_0$, we know that
$$G(s)=G(s_0)+g(s_0)s_{0}^{-(\beta-1)}\frac{1}{\beta}\left(s^\beta-s_{0}^{\beta}\right)$$
and
$$g(s)s=g(s_0)s_{0}^{-(\beta-1)}s^\beta.$$
Then, by \eqref{eq:20220102-e2}, for $s>s_0$,
\beq
g(s)s-\beta G(s)=g(s_0)s_0-\beta G(s_0)\leq 0,
\eeq
and
\begin{align}
g(s)s-\alpha G(s)=&\left(1-\frac{\alpha}{\beta}\right)g(s_0)s_{0}^{-(\beta-1)}s^\beta
+\frac{\alpha}{\beta}g(s_0)s_0-\alpha G(s_0)\nonumber\\
\geq &g(s_0)s_0-\alpha G(s_0)\geq 0.
\end{align}
Hence, $(\tilde{G}3)$ also holds.
\ep

\br\lab{remark:20210915-r1}
Let $\tilde{g}(s)$ be defined in Lemma \ref{lemma:20210915-l1}. Then:
\begin{itemize}
\item[(i)] There exists $C>0$ such that
\beq\lab{eq:20210915-e0}
\tilde{g}(s)\leq C \left(s^{1+\sigma}+s^{\beta-1}\right), \forall s\in \R^+.
\eeq
\item[(ii)]For any $M>0$, there are positive constants $C_{M,1}\geq C_{M,2}$ such that
\beq\lab{eq:20210915-e00}
C_{M,2} s^{1+\sigma}\leq \tilde{g}(s)\leq C_{M,1} s^{1+\sigma}, \forall s\in [0,M].
\eeq
\end{itemize}
\er

Fix $\alpha, \beta$   and  let $\tilde{g}$ be given in Lemma \ref{lemma:20210915-l1}. In what follows,  we focus on discussing the equation \eqref{eq:20210915-e1}.

\br
\begin{description}
\item[(i)]
For $\lambda>0$, if $u\geq 0$ is a classical solution to \eqref{eq:20210915-e1}, then the maximum principle shows that $u>0$ in $\R^N$. Furthermore, if $u(x)\rightarrow 0$ as $|x|\rightarrow \infty$, then it is easy to know that $u(x)$ decays exponentially at $\infty$ and so that $u\in H^1(\R^N)$. Conversely, if $u\in H^1(\R^N)$ is a weak solution to \eqref{eq:20210915-e1}, under the assumptions (G1),(G2) and (Q), combining with Lemma \ref{lemma:20210915-l1}, the standard elliptic estimate is applicable, and then $u$ is also a classical solution.
\item[(ii)]For $\lambda=0$, if $0\leq u\in H^1(\R^N)$ is a weak solution to \eqref{eq:20210915-e1}, under the assumptions (G1),(G2) and (Q), we still have that $u$ is a classical solution. However, if $u\geq 0$ is a classical solution to \eqref{eq:20210915-e1}, whether $u\in H^1(\R^N)$ is of interest in itself. In this respect, we have the following  Lemma \ref{lemma:20210924-l1}.
\end{description}
\er

\bl\lab{lemma:20210924-l2}(\cite[Lemma 2.2]{ArmstrongSirakov2011}) If $h(x)\in L^\infty(B_3\backslash B_{\frac{1}{2}})$ is nonnegative, and $u\geq 0$ satisfies
$$-\Delta u\geq h(x)~\hbox{in}~B_3\backslash B_{\frac{1}{2}},$$
then, there exists a constant $\bar{c}>0$ depending only on $N$ such that, for each $\Omega\subset B_2\backslash B_{1}$,
$$\inf_{B_2\backslash B_{1}}u\geq \bar{c}|\Omega| \inf_{\Omega}h.$$
\el

\bl\lab{lemma:20210924-l1}
Assume that (G1) and (G2) hold and $C_1\leq q(x)\leq C_2$, for some $C_1,C_2>0$ and any $x\in\R^N$. If $u$ is a classical solution of
\beq\lab{eq:20210924-xe1}
\begin{cases}
-\Delta u=q(x)\tilde{g}(u)~\hbox{in}~\R^N,\\
0\leq u(x)\leq M, \forall x\in \R^N,\\
u\in C^2(\R^N),
\end{cases}
\eeq
then $u\in H^1(\R^N)$, where $\tilde{g}$ is given in Lemma \ref{lemma:20210915-l1}.
\el
\bp
Suppose that $u\not\equiv 0$, then $u$ is positive in $\R^N$ by the maximum principle.
Since $u$ is bounded, by Remark \ref{remark:20210915-r1}-(ii), there exists $C>0$ such that
\beq\lab{eq:20210924-xe2}
-\Delta u=q(x)\tilde{g}(u)\geq  C u^{1+\sigma}.
\eeq
Set $u_r(x):=u(r x)$ and $r>0$, then
\beq\lab{eq:20210924-xe3}
-\Delta u_r(x)=r^2 (-\Delta u)(rx)\geq Cr^2 u(rx)^{1+\sigma}=:h(x).
\eeq
Taking $\Omega=B_2\backslash B_1$, it follows from Lemma \ref{lemma:20210924-l2} that
\beq\lab{eq:20210924-xe4}
\inf_{B_{2r}\backslash B_r} u\geq C r^2\left(\inf_{B_{2r}\backslash B_r} u\right)^{1+\sigma},
\eeq
where $C$ depends only on $M$ and $N$.
Hence,
\beq\lab{eq:20210924-xe5}
\inf_{B_{2r}\backslash B_r} u \leq C r^{-\frac{2}{\sigma}}.
\eeq

On the other hand, let $d(x):=-q(x)\frac{\tilde{g}(u(x))}{u(x)}\in L^\infty(\R^N)$ and
$$L=-\Delta+d(x),$$
then $Lu=0, u\geq 0$ in $\R^N$.
By \cite[Theorem 8.20]{GilbargTrudinger1998}, there exists some $C$ only depending on $N$ such that
\beq\lab{eq:20210924-xe6}
\sup_{B_{2r}\backslash B_r} u \leq C \inf_{B_{2r}\backslash B_r} u.
\eeq
Then, for some $C>0$ depending only on $M,N$,
\beq\lab{eq:20210924-xe7}
 u(x)\leq C |x|^{-\frac{2}{\sigma}}, \forall |x|\geq 1.
\eeq
Then, a direct computation shows that
\begin{align*}
\int_{|x|\geq 1} u^2 dx<+\infty,
\end{align*}
here $\sigma<\frac{4}{N}$ is required.
So $u\in L^2(\R^N)$ and also $u\in L^{2+\sigma}(\R^N)$, which implies that $q(\cdot)\tilde{g}(u)u\in L^1(\R^N)$. Hence, $|\nabla u|\in L^2(\R^N)$ and thus $u\in H^1(\R^N)$.
\ep

\bl\lab{lemma:20210915-l2}
Assume ($\tilde{G}1$)-($\tilde{G}3$) and (Q) hold.
 If $u\in H^1(\R^N)$ is a nonnegative solution of \eqref{eq:20210915-e1} with $\lambda\in (0,M]$, then there exists $C_M>0$, depending only on $M$, such that
$$\|u\|_\infty \leq C_M.$$
\el
\bp
We argue by contradiction. If there exists a sequence $\{\lambda_n\}\subset (0,M]$ with $\{u_n\}\subset H^1(\R^N)$, and $u_n\geq 0$ solves \eqref{eq:20210915-e1} with $\lambda=\lambda_n$, such that
$$M_n:=\|u_n\|_\infty\rightarrow +\infty.$$
Let $x_n\in \R^N$ be such that $u_n(x_n)=M_n$. Now we perform a rescaling, setting $x-x_n=\frac{y}{M_{n}^{\frac{\beta-2}{2}}}$ and
$$\tilde{u}_n(y):=\frac{1}{M_n}u_n\left(x_n+\frac{y}{M_{n}^{\frac{\beta-2}{2}}}\right),$$
then
$\displaystyle \tilde{u}_n(0)=\max_{y\in \R^N}\tilde{u}_n(y)=1$
and
\beq\lab{eq:20210915-e2}
-\Delta \tilde{u}_n(y)=q\left(\frac{y}{M_{n}^{\frac{\beta-2}{2}}}+x_n\right)\frac{\tilde{g}(M_n\tilde{u}_n(y))}{M_{n}^{\beta-1}}
-\frac{\lambda_n}{M_{n}^{\beta-2}} \tilde{u}_n(y).
\eeq
Up to a subsequence, we may assume that
\beq\lab{eq:20210915-e3}
q(x_n)\rightarrow \tau\in [q_\infty, q_0]~\hbox{as}~n\rightarrow +\infty.
\eeq
By Remark \ref{remark:20210915-r1}-(i), under the assumption (Q) and $\lambda_n\leq M, M_n\rightarrow +\infty$, the right hand side of \eqref{eq:20210915-e2} is of $L^\infty(\R^N)$. Then a standard elliptic estimate gives $\tilde{u}_n\rightarrow \tilde{u}$ in $C_{loc}^{2}(\R^N)$, where $\tilde{u}(0)=1=\max_{y\in \R^N}\tilde{u}(y)$ and
\beq\lab{eq:20210915-e4}
-\Delta \tilde{u}=\tau B \tilde{u}^{\beta-1}\;\hbox{in}~\R^N.
\eeq
So, by \cite[Theorem 8.4]{QuittnerSouplet2007} for $N=1,2$ and by \cite[Theorem 2-(ii)]{ChenLi1991} for $N\geq 3$  (see also \cite[Theorem 1.1]{GidasSpruck1981}), we know that $\tilde{u}\equiv 0$, which is a contradiction to $\tilde{u}(0)=1$.
\ep

\bl\lab{lemma:20210915-l3}
Assume ($\tilde{G}1$)-($\tilde{G}3$) and (Q) hold and let $\{u_n\}_{n=1}^{\infty}\subset H^1(\R^N)$ be positive solutions to \eqref{eq:20210915-e1} with $\lambda=\lambda_n\rightarrow 0^+$, then
\beq\lab{eq:20210915-e5}
\liminf_{n\rightarrow +\infty}\frac{\|u_n\|_{\infty}^{\sigma}}{\lambda_n}>0.
\eeq
Furthermore, if one of the following conditions holds
\begin{itemize}
\item[(i)]$N\leq 4$;
\item[(ii)]$N\geq 5$ and $\sigma\leq \frac{2}{N-2}$;
\item[(iii)] $N\geq 5, \sigma\in (\frac{2}{N-2}, \frac{4}{N})$ and $\|u_n\|_{L^2(\R^N)}\stackrel{n}{\rightarrow} 0$;
\end{itemize}
then
\beq\lab{eq:20210915-e6}
\|u_n\|_\infty\rightarrow 0~\hbox{as}~n\rightarrow +\infty.
\eeq
\el
\bp
Without loss of generality, we may assume that $\lambda_n\leq 1, \forall n\in \NN$. By Lemma \ref{lemma:20210915-l2}, there exists some $M>0$ such that
\beq\lab{eq:20210915-e7}
\|u_n\|_\infty\leq M, \forall n\in \NN.
\eeq
Let $x_n\in \R^N$ be such that $u_n(x_n)=\|u_n\|_\infty$. Define
$$\bar{u}_n(x):=\frac{1}{\|u_n\|_\infty}u_n\left(\frac{x}{\sqrt{\lambda_n}}+x_n\right),$$
then $\displaystyle 1=\bar{u}_n(0)=\max_{x\in \R^N}\bar{u}_n(x)$ and
\beq\lab{eq:20210915-e8}
-\Delta \bar{u}_n(x)+\bar{u}_n(x)=q\left(\frac{x}{\sqrt{\lambda_n}}+x_n\right)\frac{\tilde{g}(\|u_n\|_\infty \bar{u}_n(x))}{\lambda_n\|u_n\|_\infty}.
\eeq
Taking $x=0$, by the maximum principle and Remark \ref{remark:20210915-r1}-(ii), we have that
\begin{align*}
1=\bar{u}_n(0)\leq& -\Delta \bar{u}_n(0)+\bar{u}_n(0)\\
=&q(x_n)\frac{\tilde{g}(\|u_n\|_\infty )}{\lambda_n\|u_n\|_\infty}
\leq q_0 C_M  \frac{\|u_n\|_\infty^\sigma}{\lambda_n}.
\end{align*}
Hence,
$$\liminf_{n\rightarrow +\infty}\frac{\|u_n\|_\infty^\sigma}{\lambda_n}\geq \frac{1}{q_0 C_M}>0.$$
To prove \eqref{eq:20210915-e6}, we argue by contradiction for the cases (i) and (ii). Suppose that $\displaystyle\liminf_{n\rightarrow +\infty}\|u_n\|_\infty>0$, and let $v_n(x)=u_n(x+x_n)$. Then,
\beq\lab{eq:20210915-e9}
-\Delta v_n=q(x+x_n)\tilde{g}(v_n)-\lambda_n v_n~\hbox{in}~\R^N.
\eeq
Since $q(x+x_n)\tilde{g}(v_n)-\lambda_n v_n\in L^\infty(\R^N)$, it follows from a standard elliptic estimate that
 \[
 v_n\rightarrow v \text{ in } C_{loc}^{2}(\R^N)\text{ with } \displaystyle v(0)=\lim_{n\rightarrow \infty} u_n(x_n)>0.
 \]
Furthermore, $v$ is bounded and satisfies
\beq\lab{eq:20210915-e10}
-\Delta v\geq q_\infty \tilde{g}(v)\geq 0,\,x\in\R^N.
\eeq
Now, we claim that there always is a contradiction in both cases of (i)-(ii).

Case (i): If $N\leq 4$, we note that $\displaystyle\sigma+1<1+\frac{4}{N}$, by
Remark \ref{remark:20210915-r1}-(ii) and \eqref{eq:20210915-e10}, there exists some $C>0$ such that
\beq\lab{eq:20210924-e1}
-\Delta v\geq C v^{1+\sigma}~\hbox{in}~\R^N,
\eeq
where
$$
1+\sigma< \frac{N}{(N-2)_+}:
=\begin{cases}
\infty,\;\hbox{if}~N\leq 2,\\
\frac{N}{N-2}\;~\hbox{if}~N\geq 3.
\end{cases}
$$
Hence, by \cite[Theorem 8.4]{QuittnerSouplet2007}, $v\equiv 0$ in $\R^N$ , which is a contradiction to $v(0)>0$.

Case (ii): If $N\geq 5$ and $\sigma\leq \frac{2}{N-2}$, a similar argument to case (i) implies that $v\equiv 0$ in $\R^N$, which is also a contradiction.

For case (iii), we first prove that for any $p>\max\left\{2,\frac{N}{2}\right\}$, there exists some $C(p)>0$ such that
\beq\lab{eq:20220105-be1}
\|u_n\|_\infty \leq C(p) \|u_n\|_{2}^{\frac{2}{p}}, \forall n\in \N.
\eeq
Indeed, for all $n\in \N$,
$$\Delta u_n=\lambda_n u_n -q(x)\tilde{g}(u_n).$$
By (Q), $\lambda_n\leq 1,\|u_n\|_\infty\leq M$ and \eqref{eq:20210915-e0}, there exists $K>0$ such that
\beq\lab{eq:20220105-be2}
|\Delta u_n|\leq K|u_n|.
\eeq
Using the Calderon-Zymund estimate \cite[Chap. 2,3, Prop.8]{DautrayLions1984}:
\beq\lab{eq:20220105-be3}
\left\|\frac{\partial^2 u}{\partial x_i\partial x_j}\right\|_p \leq C(p) \|\Delta u\|_p~\hbox{for all}~u\in W^{2,p}(\R^N), 1<p<\infty,
\eeq
we see from \eqref{eq:20220105-be2} that if $\{u_n\}\subset L^p(\R^N)$ for some $1<p<\infty$, then it also holds that $\{u_n\}\subset W^{2,p}(\R^N)$ and
\beq\lab{eq:20220105-be4}
\|u_n\|_{W^{2,p}(\R^N)}\leq C(p)\|u_n\|_p, ~\forall n\in \N.
\eeq
Take $p>\max\{2,\frac{N}{2}\}$, by $\|u_n\|_\infty\leq M$, we obtain that
\beq\lab{eq:20220105-be5}
\|u_n\|_p\leq M^{\frac{p-2}{p}} \|u_n\|_{2}^{\frac{2}{p}}.
\eeq
We recall the continuous embedding $W^{2,p}(\R^N)\hookrightarrow L^\infty(\R^N)$ for any $p>\frac{N}{2}$. Hence,
\beq\lab{eq:20220105-be6}
\|u_n\|_\infty \leq C_1(p)\|u_n\|_{W^{2,p}(\R^N)} \leq C_2(p) \|u_n\|_p
\leq C_3(p) \|u_n\|_{2}^{\frac{2}{p}},
\eeq
we prove \eqref{eq:20220105-be1}.

Hence, in Case (iii), that is,  $N\geq 5$, $\frac{2}{N-2}<\sigma<\frac{4}{N}$ and $\|u_n\|_2\stackrel{n}{\rightarrow} 0$,  it follows \eqref{eq:20220105-be1} that $\|u_n\|_\infty \stackrel{n}{\rightarrow} 0$ also holds.
\ep

\br\lab{remark:20220105-br1}
Under the assumptions $(\tilde{G}1)-(\tilde{G}3)$ and (Q), if $\{u_n\}_{n=1}^{\infty}\subset H^1(\R^N)$ be positive solutions to \eqref{eq:20210915-e1} with $\lambda=\lambda_n\rightarrow 0^+$ such that $\|u_n\|_2\stackrel{n}{\rightarrow}0$, then $\|u_n\|_\infty\stackrel{n}{\rightarrow}0$ by Lemma \ref{lemma:20210915-l3}.
\er

For $\lambda>0$ and $u\in H^1(\R^N)$, define
\beq\lab{eq:20210915-e12}
\tilde{I}_{\infty,\lambda}[u]:=\frac{1}{2}\|\nabla u\|_2^2+\lambda \|u\|_2^2-q_\infty\int_{\R^N}\tilde{G}(u)dx
\eeq
and the associate least energy
\beq\lab{eq:20210915-e13}
\ell_{\infty,\lambda}:=\inf\left\{\tilde{I}_{\infty,\lambda}[u]: u\in H^1(\R^N)\backslash \{0\}~\hbox{with}~\tilde{I}'_{\infty,\lambda}[u]=0\right\},
\eeq
the mountain pass level
\beq\lab{eq:20210901-e5}
m_{\infty,\lambda}:=\inf_{\gamma\in \Gamma_\infty}\max_{t\in [0,1]}\tilde{I}_{\infty,\lambda}(\gamma(t)),
\eeq
where
\beq\lab{eq:20210901-e6}
\Gamma_\infty:=\left\{\gamma\in C([0,1], H^1(\R^N)): \gamma(0)=0, \tilde{I}_{\infty,\lambda}(\gamma(1))<0.\right\}.
\eeq

Then, we have

\bl\lab{lemma:20210915-bl3}
If ($\tilde{G}1$)-($\tilde{G}3$) and (Q) hold,
then  $\ell_{\infty,\lambda}>0$ for any $\lambda>0$. In particular, $\ell_{\infty,\lambda}=m_{\infty,\lambda}$ is nondecreasing with respect to $\lambda\in \R^+$.
\el
\bp
We mention that the existence of a least energy solution for
\beq\lab{eq:20220102-we1}
-\Delta u+u=q_\infty \tilde{g}(u)~\hbox{in}~\R^N, u\in H^1(\R^N),
\eeq
 was proven by Berestycki and Lions in \cite{BerestyckiLions1983}. It was also shown that the least energy solution is of mountain pass type (see, e.g.,\cite{JeanTanaka2002} for $N\geq 2$ and \cite{JeanTanaka2003} for $N=1$). We omit the proof here.
\ep

\s{Asymptotic behavior of mountain pass solution}\lab{sec:behavior}
\renewcommand{\theequation}{3.\arabic{equation}}
Define
\beq\lab{eq:bu-zzze1}
\tilde{I}_\lambda[u]:=\frac{1}{2}\|\nabla u\|_2^2+\frac{1}{2}\lambda \|u\|_2^2-\int_{\R^N}q(x)\tilde{G}(u)dx,\,\,u\in H^1(\R^N),
\eeq
and
\beq\lab{eq:20210901-e5}
m_\lambda:=\inf_{\gamma\in \Gamma}\max_{t\in [0,1]}\tilde{I}_\lambda(\gamma(t)),
\eeq
where $\Gamma$ is defined by
\beq\lab{eq:20210917-bue1}
\Gamma:=\left\{\gamma\in C([0,1], H^1(\R^N)): \gamma(0)=0, \tilde{I}_\lambda(\gamma(1))<0\right\}.
\eeq

\bl\lab{lemma:nb-l1}
Let ($\tilde{G}1$)-($\tilde{G}3$) and (Q) be satisfied. Then,
for any $\lambda>0$,
$$0<m_\lambda<m_{\infty,\lambda}$$
and $m_{\infty,\lambda}\rightarrow 0$ as $\lambda\rightarrow 0^+$.
\el
\bp
For $\lambda>0$, $m_{\infty,\lambda}$ is attained (see the proof of Lemma \ref{lemma:20210915-bl3}). Let $u_{\infty,\lambda}$ be a mountain pass solution for \eqref{eq:20220102-we1}, and by the regularity result, we may suppose that $u_{\infty,\lambda}>0$ in $\R^N$. And there exists some $\gamma\in \Gamma_\infty$ such that
$\gamma(t_0)=u_{\infty,\lambda}$ for some $t_0\in (0,1]$ with
\beq
\max_{t\in [0,1]}\tilde{I}_{\infty,\lambda}(\gamma(t))=\tilde{I}_{\infty,\lambda}(\gamma(t_0))
=\tilde{I}_{\infty,\lambda}(u_{\infty,\lambda})=m_{\infty,\lambda}.
\eeq
We note that $\gamma\in \Gamma$ since $\tilde{I}_\lambda(\gamma(1))\leq \tilde{I}_{\infty,\lambda}(\gamma(1))<0$.  Then there exists $t_1\in (0,1)$ such that
\begin{align*}
m_\lambda\leq& \max_{t\in [0,1]}\tilde{I}_{\lambda}(\gamma(t))
=\tilde{I}_{\lambda}(\gamma(t_1))
<\tilde{I}_{\infty,\lambda}(\gamma(t_1))\\
\leq&\max_{t\in [0,1]}\tilde{I}_{\infty,\lambda}(\gamma(t))
=m_{\infty,\lambda}.
\end{align*}
On the other hand, for $\lambda>0$ fixed, the norm $\|u\|_\lambda$ defined by
$$\|u\|_\lambda:=\left(\|\nabla u\|_2^2+\lambda\|u\|_2^2\right)^{\frac{1}{2}},$$
is equivalent to the usual norm $\|u\|_{H^1(\R^N)}=\|u\|$.
By the Sobolev embeddings, for $\|u\|_\lambda=\rho>0$ small, we have
\begin{align*}
\tilde{I}_\lambda[u]=&\frac{1}{2}\|u\|_\lambda^2-\int_{\R^N}q(x)\tilde{G}(u)dx\\
\geq &\frac{1}{2}\|u\|_\lambda^2- C \left(\|u\|_\alpha^\alpha+\|u\|_\beta^\beta\right)\\
\geq& \frac{1}{2}\|u\|_\lambda^2 -C \left(\frac{1}{2}\|u\|_\lambda^\alpha +\frac{1}{2}\|u\|_\lambda^\beta\right).
\end{align*}
Since $\alpha,\beta>2$, it is easy to see
$$m_\lambda\geq \inf_{\|u\|_\lambda=\rho}\tilde{I}_\lambda[u]>0.$$

For $\lambda>0$ small, let $u_{\infty,\lambda}$ be the corresponding mountain pass solution of \eqref{eq:20220102-we1}. Under the conditions $(\tilde{G}1)$ and $(\tilde{G}2)$, we remark that the assumption (G3) required in \cite{JeanZhangZhong2021} holds automatically, see \cite[Remark 2.6]{JeanZhangZhong2021}. So by \cite[Theorem 5.1]{JeanZhangZhong2021}, $u_{\infty,\lambda}$ is the unique nontrivial nonnegative solution. Furthermore, by \cite[Theorem 1.3]{JeanZhangZhong2021}, we have that $u_{\infty,\lambda}\rightarrow 0$ in $H^1(\R^N)$. Hence,
$$m_{\infty,\lambda}=\tilde{I}_{\infty,\lambda}[u_{\infty,\lambda}]\rightarrow 0~\hbox{as}~\lambda\rightarrow 0^+.$$
\ep

\bl\lab{lemma:nb-l2} Assume ($\tilde{G}1$)-($\tilde{G}3$) and (Q) hold. For $m_{\infty,\lambda}$ defined by \eqref{eq:20210901-e5}, there holds
$$\displaystyle \lim_{\lambda\rightarrow 0^+}\frac{m_{\infty,\lambda}}{\lambda}=0.$$
\el
\bp
Noting that
\beq\lab{eq:20210916-e1}
m_{\infty,\lambda}=\frac{1}{2}\|\nabla u_{\infty,\lambda}\|_2^2+\frac{1}{2}\lambda \|u_{\infty,\lambda}\|_2^2-q_\infty \int_{\R^N}\tilde{G}(u_{\infty,\lambda})dx,
\eeq

\beq\lab{eq:20210916-e2}
\|\nabla u_{\infty,\lambda}\|_2^2+\lambda \|u_{\infty,\lambda}\|_2^2=q_\infty \int_{\R^N}\tilde{g}(u_{\infty,\lambda})u_{\infty,\lambda} dx,
\eeq

\beq\lab{eq:20210916-e3}
\|\nabla u_{\infty,\lambda}\|_2^2=N q_\infty\int_{\R^N}\left[\frac{1}{2}\tilde{g}(u_{\infty,\lambda})u_{\infty,\lambda}
-\tilde{G}(u_{\infty,\lambda})\right]dx
\eeq
and
\beq\lab{eq:20210916-e4}
\alpha \int_{\R^N}\tilde{G}(u_{\infty,\lambda})dx\leq\int_{\R^N}\tilde{g}(u_{\infty,\lambda})u_{\infty,\lambda} dx\leq \beta \int_{\R^N}\tilde{G}(u_{\infty,\lambda})dx.
\eeq
Then,
\beq\lab{eq:20210916-e5}
\begin{cases}
Nq_\infty \frac{\alpha-2}{2}\int_{\R^N}\tilde{G}(u_{\infty,\lambda})dx
\leq \|\nabla u_{\infty,\lambda}\|_2^2\leq Nq_\infty \frac{\beta-2}{2}\int_{\R^N}\tilde{G}(u_{\infty,\lambda})dx,\\
q_\infty \frac{\alpha-2}{2}\int_{\R^N}\tilde{G}(u_{\infty,\lambda})dx
\leq m_{\infty,\lambda}\leq q_\infty \frac{\beta-2}{2}\int_{\R^N}\tilde{G}(u_{\infty,\lambda})dx,\\
q_\infty\left[\alpha-\frac{N(\beta-2)}{2}\right]\int_{\R^N}\tilde{G}(u_{\infty,\lambda})dx
\leq \lambda \|u_{\infty,\lambda}\|_2^2\leq
q_\infty\left[\beta-\frac{N(\alpha-2)}{2}\right]\int_{\R^N}\tilde{G}(u_{\infty,\lambda})dx.
\end{cases}
\eeq
By $2<\alpha\leq \beta<2+\frac{4}{N}$, we know that
$$0<\alpha-\frac{N(\beta-2)}{2}\leq \beta-\frac{N(\alpha-2)}{2}.$$
Hence, by \cite[Theorem 1.3]{JeanZhangZhong2021} again, $u_{\infty,\lambda}\rightarrow 0$ in $H^1(\R^N)$ as $\lambda\rightarrow 0^+$. Then it follows \eqref{eq:20210916-e5} that
\begin{align*}
\frac{m_{\infty,\lambda}}{\lambda}\leq& \frac{\beta-2}{2}q_\infty \frac{\int_{\R^N}\tilde{G}(u_{\infty,\lambda})dx}{\lambda}
\leq \frac{\frac{\beta-2}{2}q_\infty}{q_\infty\left[\alpha-\frac{N(\beta-2)}{2}\right]} \frac{\lambda \|u_{\infty,\lambda}\|_2^2}{\lambda}\\
=&\frac{\beta-2}{N\left(2+\frac{2}{N}\alpha -\beta\right)}\|u_{\infty,\lambda}\|_2^2\rightarrow 0~\hbox{as}~\lambda\rightarrow 0^+.
\end{align*}
\ep

\bc\lab{cro:nb-c1}
If ($\tilde{G}1$)-($\tilde{G}3$) and (Q) hold, then, for $m_\lambda$ is given by \eqref{eq:20210901-e5}, we have
$$\displaystyle \lim_{\lambda\rightarrow 0^+}\frac{m_{\lambda}}{\lambda}=0.$$
\ec
\bp
It follows from Lemmas \ref{lemma:nb-l1} and \ref{lemma:nb-l2}.
\ep

\bo\lab{prop:nb-p2}
Under the conditions ($\tilde{G}1$)-($\tilde{G}3$) and (Q),
for $\lambda>0$, let $u_\lambda$ be a mountain pass solution of \eqref{eq:20210915-e1}, then
$$u_\lambda\rightarrow 0~\hbox{in}~H^1(\R^N)~\hbox{as}~\lambda\rightarrow 0^+.$$
Consequently, $\lambda=0$ is a bifurcation point for \eqref{eq:20210915-e1}.
\eo
\bp
Since $u_\lambda$ is a solution to \eqref{eq:20210915-e1},
\beq\lab{eq:20210916-e6}
\|\nabla u_\lambda\|_2^2+\lambda \|u_\lambda\|_2^2
=\int_{\R^N}q(x)\tilde{g}(u_\lambda)u_\lambda dx.
\eeq
Hence,
\begin{align*}
m_\lambda=&\frac{1}{2}\left(\|\nabla u_\lambda\|_2^2+\lambda \|u_\lambda\|_2^2\right)-\int_{\R^N}q(x)\tilde{G}(u_\lambda)dx\\
\geq &\frac{1}{2}\left(\|\nabla u_\lambda\|_2^2+\lambda \|u_\lambda\|_2^2\right)-\frac{1}{\alpha}\int_{\R^N}q(x)\tilde{g}(u_\lambda)u_\lambda dx\\
=&\left(\frac{1}{2}-\frac{1}{\alpha}\right)\left(\|\nabla u_\lambda\|_2^2+\lambda \|u_\lambda\|_2^2\right).
\end{align*}
So, for $\lambda\leq 1$, we have
\beq\lab{eq:20210916-e7}
m_\lambda\geq \left(\frac{1}{2}-\frac{1}{\alpha}\right) \lambda \|u_{\lambda}\|_{H^1(\R^N)}^{2}.
\eeq
By Corollary \ref{cro:nb-c1}, we  obtain that $$\lim_{\lambda\rightarrow 0^+}\|u_{\lambda}\|_{H^1(\R^N)}^{2}\leq \frac{2\alpha}{\alpha-2}\lim_{\lambda\rightarrow 0^+}\frac{m_\lambda}{\lambda}=0.$$
\ep

\s{Existence of mountain pass solution}\lab{sec:existence}
\renewcommand{\theequation}{4.\arabic{equation}}

Under the assumptions of ($\tilde{G}1$)-($\tilde{G}3$) and (Q), it is easy to check that $\tilde{I}_\lambda$ satisfies the mountain pass geometry. Then there exists a $(PS)_{m_\lambda}$ sequence $\{u_n\}\subset H^1(\R^N)$, i.e., as $n\rightarrow +\infty$,
\beq\lab{eq:20210916-e8}
\tilde{I}'_\lambda[u_n]\rightarrow 0, \tilde{I}_\lambda[u_n]\rightarrow m_\lambda.
\eeq
Since $\tilde{g}$ satisfies the so-called Ambrosetti-Rabinowitz type condition (see Lemma \ref{lemma:20210915-l1}-($\tilde{G}3$)), it is standard to check that $\{u_n\}$ is bounded in $H^1(\R^N)$.
Up to a subsequence, we assume that $u_n\stackrel{n}{\rightharpoonup} u$ weakly in $H^1(\R^N)$ and $u$ is a weak solution of \eqref{eq:20210915-e1}, and our following Lemma shows that $u \not\equiv 0$.

\bl\lab{lemma:20211020-wl2} Under the assumptions of ($\tilde{G}1$)-($\tilde{G}3$) and (Q), let $u$ be the weak limit of the sequence $\{u_n\}$ given by \eqref{eq:20210916-e8}, then $u\not\equiv 0$.
\el
\bp
By contradiction, if $u\equiv 0$, then
\beq\lab{eq:20211020-wbue1}
u_n\stackrel{n}{\rightarrow} 0~\hbox{ strongly in}~L_{loc}^{p}(\R^N), \forall p\in [1,2^*).
\eeq

\vskip0.02in
\noindent
{\bf Step 1.} We firstly prove that there exist $\eta>0, R>0$ and $\{y_n\}\subset \R^N$ such that
\beq\lab{eq:20211020-we2}
\liminf_{n\rightarrow +\infty}\int_{B_R(y_n)} |u_n|^2 dx\geq \eta.
\eeq
 Indeed, for any $R>0$, if
$$\sup_{y\in \R^N} \int_{B_R(y)}|u_n|^2 dx\rightarrow 0\;\hbox{as}\;n\rightarrow \infty.$$
Then, we infer that $u_{n}\rightarrow 0$ strongly in $L^p(\R^N)$ for any $2<p<2^*$ (See \cite[Lemma 1.21]{Willem1996}).
Then, by Lemma \ref{lemma:20210915-l1}-($\tilde{G}3$) we have
\beq\lab{eq:20210916-e9}
\lim_{n\rightarrow \infty}\int_{\R^N}q(x)\tilde{G}(u_n)dx=0=\lim_{n\rightarrow \infty}\int_{\R^N}q(x)\tilde{g}(u_n)u_n dx.
\eeq
These imply that $u_n \stackrel{n}\rightarrow 0$ strongly in $H^1(\R^N)$ and hence
$m_\lambda=\lim_{n\rightarrow \infty}\tilde{I}_\lambda[u_n]=0$, which contradicts Lemma \ref{lemma:nb-l1}.

\vskip0.02in
\noindent
{\bf Step 2.} Let $\{y_n\}$ be obtained in Step 1, then $\{y_n\}$ is unbounded in $\R^N$.\\
If there exists $M>0$ such that $\sup_{n}|y_n|\leq M$, then, for $R>0$ given in Step 1, we have
\beq\lab{eq:20211020-wbue2}
\int_{B_{R+M}(0)}|u_n|^2 dx\geq \int_{B_R(y_n)}|u_n|^2 dx\geq \eta>0,
\eeq
which is a contradiction to our assumption \eqref{eq:20211020-wbue1}.

\vskip0.02in
\noindent
{\bf Step 3.} We claim that $(q(x)-q_\infty)\tilde{g}(u_n)\rightarrow 0$ in $H^{-1}$ as $n\rightarrow\infty$.\\
Indeed, for any $\varphi\in H^1(\R^N)$, it follows from \eqref{eq:20210915-e0}, H\"older inequality and Sobolev embeddings that, for any $R>0$,
\begin{align*}
&\left|\int_{\R^N}(q(x)-q_\infty)\tilde{g}(u_n) \varphi dx\right|\\
\leq&\int_{|x|\leq R} |q(x)-q_\infty| \tilde{g}(u_n) |\varphi| dx +\int_{|x|\geq R} |q(x)-q_\infty| \tilde{g}(u_n) |\varphi| dx\\
\leq&C \left[\left(\int_{|x|\leq R} |u_n|^{2+\sigma}dx\right)^{\frac{1+\sigma}{2+\sigma}} \|\varphi\|_{2+\sigma} + \left(\int_{|x|\leq R} |u_n|^\beta dx\right)^{\frac{\beta-1}{\beta}} \|\varphi\|_\beta \right] \\
&+\sup_{|x|\geq R} |q(x)-q_\infty| \int_{|x|\geq R} C(|u_n|^{1+\sigma} +|u_n|^\beta) |\varphi| dx\\
=&o(1)\|\varphi\|_{H^1},~\hbox{by letting $n\rightarrow +\infty$ and using \eqref{eq:20211020-wbue1}, then letting $R\rightarrow +\infty$}.
\end{align*}

\vskip0.02in
\noindent
{\bf Step 4.}
Put $\omega_n(x):=u_n(x+y_n)$, then $\tilde{I}'_{\infty,\lambda}[\omega_n]\rightarrow 0$ in $H^{-1}$.\\
For any $\varphi\in H^1(\R^N)$, denote $\varphi_n(x)=\varphi(x-y_n)$. Since $\|\varphi_n\|_{H^1}=\|\varphi\|_{H^1}$,
\begin{align*}
\tilde{I}'_{\infty,\lambda}[\omega_n]\varphi=&\tilde{I}'_{\infty,\lambda}[u_n]\varphi_n
=\tilde{I}'_\lambda[u_n]\varphi_n+ \int_{\R^N}(q(x)-q_\infty) \tilde{g}(u_n)\varphi_n dx.
\end{align*}
By Step 3, $(q(x)-q_\infty)\tilde{g}(u_n)\rightarrow 0$ in $H^{-1}$, then it follows from $\tilde{I}'_\lambda[u_n]\rightarrow 0$ in $H^{-1}$ that
$$|\tilde{I}'_{\infty,\lambda}[\omega_n]\varphi|\leq o(1)\|\varphi\|_{H^1}.$$
That is, $\tilde{I}'_{\infty,\lambda}[\omega_n]\rightarrow 0$ in $H^{-1}$.

\vskip0.02in
\noindent
{\bf Step 5.} Up to a subsequence, $\omega_n\stackrel{n}{\rightharpoonup} \omega$ in $H^1(\R^N)$. By Step 1 and Step 4, we know that $\omega\not\equiv 0$ is a weak solution to equation:
\beq\lab{eq:20211020-we8}
-\Delta \omega+\lambda \omega=q_\infty \tilde{g}(\omega)~\hbox{in}~\R^N.
\eeq
By \eqref{eq:20210915-e13} and Lemma \ref{lemma:20210915-bl3}, we have
\beq\lab{eq:20211020-wbue3}
\tilde{I}_{\infty,\lambda}[\omega]\geq \ell_{\infty,\lambda}=m_{\infty,\lambda}.
\eeq

\vskip0.02in
\noindent
{\bf Step 6.} By the well known Brezis-Lieb lemma, we have
\begin{align*}
\tilde{I}'_{\infty,\lambda}[\omega_n-\omega] (\omega_n-\omega)=&\|\omega_n-\omega\|_\lambda^2-q_\infty\int_{\R^N}\tilde{g}(\omega_n-\omega) (\omega_n-\omega)dx\\
=&\|\omega_n\|_\lambda^2-\|\omega\|_\lambda^2-q_\infty \int_{\R^N}[\tilde{g}(\omega_n)\omega_n -\tilde{g}(\omega)\omega]dx+o(1)\\
=&\tilde{I}'_{\infty,\lambda}[\omega_n]\omega_n-\tilde{I}'_{\infty,\lambda}[\omega]\omega +o(1)\\
=&\tilde{I}'_{\infty,\lambda}[\omega_n]\omega_n+o(1)\\
=&o(1),
\end{align*}
where the fact that $\{\omega_n\}$ is bounded in $H^1(\R^N)$ and the result of Step 4 are used.
Then, by $(\tilde{G}3)$,
\begin{align*}
\tilde{I}_{\infty,\lambda}[\omega_n-\omega]=&\frac{1}{2}\|\omega_n-\omega\|_\lambda^2-\int_{\R^N}q_\infty \tilde{G}(\omega_n-\omega)dx\\
\geq &\frac{1}{2}\|\omega_n-\omega\|_\lambda^2-\frac{1}{\alpha}\int_{\R^N}q_\infty \tilde{g}(\omega_n-\omega) (\omega_n-\omega)dx\\
=&\left(\frac{1}{2}-\frac{1}{\alpha}\right)\|\omega_n-\omega\|_\lambda^2+o(1)\\
\geq&o(1).
\end{align*}
Applying the Brezis-Lieb lemma again,
\beq
\liminf_{n\rightarrow \infty}\tilde{I}_{\infty,\lambda}[\omega_n]=\liminf_{n\rightarrow \infty}\left[\tilde{I}_{\infty,\lambda}[\omega]+\tilde{I}_{\infty,\lambda}[\omega_n-\omega]+o(1)\right]
\geq\tilde{I}_{\infty,\lambda}[\omega]\geq m_{\infty,\lambda}.
\eeq
On the other hand, by Step 3 and $(\tilde{G}3)$, we have
\begin{align*}
\left|\int_{\R^N}(q(x)-q_\infty)\tilde{G}(u_n)dx\right|\leq \frac{1}{\alpha}\int_{\R^N}(q(x)-q_\infty)\tilde{g}(u_n)u_ndx=o(1)\|u_n\|_{H^1}=o(1).
\end{align*}
Hence,
\begin{align*}
m_\lambda=&\tilde{I}_\lambda[u_n]+o(1)\\
=&\frac{1}{2}\|\nabla \omega_n\|_2^2+\frac{1}{2}\lambda\|\omega_n\|_2^2-\int_{\R^N}q(x)\tilde{G}(u_n)dx+o(1)\\
=&\tilde{I}_{\infty,\lambda}[\omega_n]-\int_{\R^N}q(x)\tilde{G}(u_n)dx+\int_{\R^N}q_\infty\tilde{G}(\omega_n)dx+o(1)\\
=&\tilde{I}_{\infty,\lambda}[\omega_n]-\int_{\R^N}(q(x)-q_\infty)\tilde{G}(u_n)dx+o(1)\\
=&\tilde{I}_{\infty,\lambda}[\omega_n]+o(1)\geq\ell_{\infty,\lambda}+o(1)\\
=&m_{\infty,\lambda}+o(1),
\end{align*}
which is a contradiction to Lemma \ref{lemma:nb-l1} that $m_\lambda<m_{\infty,\lambda}$. So, Lemma \ref{lemma:20211020-wl2} is proved.
\ep
\bt\lab{thm:mountain-pass}
Let ($\tilde{G}1$)-($\tilde{G}3$) and (Q) be satisfied.
For any $\lambda>0$, there exists a mountain pass solution $0<u_\lambda\in H^1(\R^N)$ to \eqref{eq:20210915-e1}, that is,
\beq
\tilde{I}'_\lambda(u_\lambda)=0~\hbox{and}~\tilde{I}_\lambda(u_\lambda)=m_\lambda.
\eeq
\et
\bp
For $\{u_n\}\subset H^1(\R^N)$ given by \eqref{eq:20210916-e8}, we claim that $u_n\rightarrow u$ in $H^1(\R^N)$. Otherwise, for $\phi_n(x):=u_n(x)-u(x)$, we have $\phi_n\rightharpoonup 0$ in $H^1(\R^N)$ but $\phi_n\not\rightarrow 0$ in $H^1(\R^N)$. By the well known Brezis-Lieb lemma,
\beq\lab{eq:20210916-e12}
\begin{cases}
\int_{\R^N}q(x)\tilde{g}(u_n)u_n dx=\int_{\R^N}q(x)\tilde{g}(u)u dx +\int_{\R^N}q(x)\tilde{g}(\phi_n)\phi_n dx+o(1),\\
\|\nabla u_n\|_2^2+\lambda\|u_n\|_2^2=\|\nabla u\|_2^2+\lambda\|u\|_2^2+\|\nabla \phi_n\|_2^2+\lambda\|\phi_n\|_2^2+o(1),\\
\|\nabla u\|_2^2+\lambda\|u\|_2^2=\int_{\R^N}q(x)\tilde{g}(u)u dx,\\
\|\nabla u_n\|_2^2+\lambda\|u_n\|_2^2=\int_{\R^N}q(x)\tilde{g}(u_n)u_n dx+o(1).
\end{cases}
\eeq
Then, there exist $\eta_1>0,R_1>0$ and $\{z_n\}\subset \R^N$ with $|z_n|\rightarrow +\infty$ such that
\beq\lab{eq:20210916-e13}
\liminf_{n\rightarrow \infty} \int_{B_{R_1}(z_n)} |\phi_{n}|^2 dx\geq \eta_1>0.
\eeq
Otherwise, similar to the proof of \eqref{eq:20210916-e9}, we have $\displaystyle \int_{\R^N}q(x)\tilde{g}(\phi_n)\phi_n dx=o(1)$, and then \eqref{eq:20210916-e12} implies that
$u_n\rightarrow u$ in $H^1(\R^N)$, which is a contradiction.

So, by \eqref{eq:20210916-e13} and $\phi_n\rightharpoonup 0$ in $H^1(\R^N)$, we have, for some $\psi\in H^1(\R^N)$,
$$\psi_n(x):=\phi_n(x+z_n)\stackrel{n}{\rightharpoonup} \psi\not\equiv 0\,\,\mbox{in}\,\,H^1(\R^N).$$
Then by applying a similar argument in Lemma \ref{lemma:20211020-wl2} (see Step 4), we get that $\psi_n$ is a $(PS)$ sequence of $\tilde{I}_{\infty,\lambda}$, and thus
\beq
\tilde{I}_{\infty,\lambda}[\psi_n]\geq \ell_{\infty,\lambda}+o(1)=m_{\infty,\lambda}+o(1).
\eeq
Hence, by the Brezis-Lieb lemma again,
\begin{align*}
m_\lambda=&\tilde{I}_\lambda[u_n]+o(1)\\
=&\tilde{I}_\lambda[u]+\tilde{I}_\lambda[\phi_n]+o(1)\\
=&\tilde{I}_\lambda[u]+\tilde{I}_{\infty,\lambda}[\psi_n]+o(1)\\
>&\tilde{I}_{\infty,\lambda}[\psi_n]+o(1)\\
\geq&m_{\infty,\lambda}+o(1),
\end{align*}
which contradicts to Lemma \ref{lemma:nb-l1} again, where we use the fact $\tilde{I}[u]>0$, see Lemma \ref{lemma:bul1} below.

Hence, $u_n\stackrel{n}{\rightarrow} u$ in $H^1(\R^N)$ and $\tilde{I}_\lambda[u]=m_\lambda$. That is, $u$ is a mountain pass solution to \eqref{eq:20210915-e1}. The proof of Theorem \ref{thm:mountain-pass} is complete.
\ep

\bl\lab{lemma:bul1}
If ($\tilde{G}1$)-($\tilde{G}3$) and (Q) hold, then,
for any fixed $\lambda>0$, there exists $\eta_\lambda>0$ such that
$$\inf\{\tilde{I}_\lambda[u]:  u\in H^1(\R^N)\backslash\{0\} ~\hbox{and}~ \langle \tilde{I}'_\lambda[u], u\rangle=0\}\geq \eta_\lambda.$$
\el
\bp
By Sobolev embeddings, under the assumptions (Q) and $(\tilde{G}3)$, there exists $\bar{C}>0$, which depends only on $\lambda,\alpha$ and $\beta$, such that
$$\int_{\R^N}q(x)\tilde{g}(u)u \leq \bar{C} (\|u\|_\lambda^\alpha+\|u\|_\lambda^\beta),\forall u\in H^1(\R^N).$$
For any $u\not\equiv 0$ with $\langle \tilde{I}'_\lambda[u], u\rangle=0$, we know that
$$\|u\|_\lambda^2=\int_{\R^N}q(x)\tilde{g}(u)u\leq \bar{C} (\|u\|_\lambda^\alpha+\|u\|_\lambda^\beta).$$
Then, by $\beta>\alpha>2$ there exists $C=C(\lambda,\alpha,\beta)>0$ such that
$$\inf\{\|u\|_\lambda: u\in H^1(\R^N)\backslash\{0\}~\hbox{and}~ \langle \tilde{I}'_\lambda[u], u\rangle=0\}\geq C.$$
Hence,
\begin{align*}
\tilde{I}_\lambda [u]=&\frac{1}{2}\|u\|_\lambda^2-\int_{\R^N}q(x)\tilde{G}(u)
\geq\frac{1}{2}\|u\|_\lambda^2-\frac{1}{\alpha} \int_{\R^N}q(x)\tilde{g}(u)u\\
=&\left(\frac{1}{2}-\frac{1}{\alpha}\right)\|u\|_\lambda^2
\geq \left(\frac{1}{2}-\frac{1}{\alpha}\right) C=:\eta_\lambda>0.
\end{align*}
\ep

\s{Proofs of Theorem \ref{thm:main} and Theorem \ref{cro:20210917-c1}}\lab{sec:proof}
\renewcommand{\theequation}{5.\arabic{equation}}
Now, we are ready to prove our main Theorems \ref{thm:main} and \ref{cro:20210917-c1}.\\

\noindent\textit{\bf Proof of Theorem \ref{thm:main}.}\\

Under the assumptions (G1)-(G2) and (Q), fix $\alpha,\beta$ such that $2<\alpha<2+\sigma<\beta<2+\frac{4}{N}$. Let $\tilde{g}$ be the modified function given in Lemma \ref{lemma:20210915-l1} such that ($\tilde{G}1$)-($\tilde{G}3$) hold. Then by Theorem \ref{thm:mountain-pass}, for any $\lambda>0$, problem \eqref{eq:20210915-e1} possesses a mountain pass solution $u_\lambda$. By Proposition \ref{prop:nb-p2}, we know that
\beq\lab{eq:zzz-zhy}
u_\lambda\rightarrow 0~\hbox{ strongly in}~H^1(\R^N)~\hbox{as}~\lambda\rightarrow 0^+.
\eeq
Furthermore, since $\tilde{g}$ is odd, we may assume that $u_\lambda$ is nonnegative in $\R^N$. Then, by the regularity result and maximum principle, $u_\lambda$ is strictly positive.
Let $s_0$ be given by Lemma \ref{lemma:20210915-l1}, because of \eqref{eq:zzz-zhy} and Remark \ref{remark:20220105-br1}, we conclude that $\|u_\lambda\|_\infty\rightarrow 0$ as $\lambda\rightarrow 0^+$. So, we can find some $\lambda_0>0$ such that
$u_\lambda(x)\leq s_0, \forall \lambda\in (0,\lambda_0), \forall x\in \R^N$. Then, the definition of $\tilde{g}$ implies that,
$$\tilde{g}(u_\lambda(x))\equiv g(u_\lambda(x))~\hbox{in}~\R^N, \forall \lambda\in (0,\lambda_0).$$
That is, for $\lambda\in (0,\lambda_0)$, $u_\lambda$ is essential a positive solution to the original problem \eqref{eq:B-P}.
So, using \eqref{eq:zzz-zhy} again, $\lambda=0$ is a bifurcation point for equation \eqref{eq:B-P}. The proof of Theorem \ref{thm:main} is completed.\\

\noindent\textit{\bf Proof of Theorem \ref{cro:20210917-c1}.}\\

By Theorem \ref{thm:main}, there exists some $\lambda_1>0$ small enough such that
for any $\lambda\in (0,\lambda_1)$ with $\lambda_1:=\frac{\lambda_0}{V_\infty}$, there exists a mountain pass solution $u_{1,\lambda}$ to
\beq\lab{eq:20210917-xe1}
-\Delta u+\lambda V_\infty u=q(x)g(u)~\hbox{in}~\R^N.
\eeq
Furthermore,
$$u_{1,\lambda}\rightarrow 0~\hbox{in}~H^1(\R^N)~\hbox{as}~\lambda\rightarrow 0^+.$$
For any $u\in H^1(\R^N)$, set
\beq\lab{eq:20210917-xe3}
I_{1,\lambda}[u]:=\frac{1}{2}\|\nabla u\|_2^2+\frac{1}{2}\lambda\int_{\R^N}V_\infty u^2 dx-\int_{\R^N}q(x)G(u)dx
\eeq
and
\beq\lab{eq:20210917-xe5}
I_{\lambda}[u]:=\frac{1}{2}\|\nabla u\|_2^2+\frac{1}{2}\lambda\int_{\R^N}V(x) u^2 dx-\int_{\R^N}q(x)G(u)dx.
\eeq
Define
\beq\lab{eq:20210917-xe6}
m_{1,\lambda}=\inf_{\gamma\in \Gamma}\max_{t\in [0,1]}I_{1,\lambda}(\gamma(t)),
\eeq
where
\beq\lab{eq:20210917-xe7}
\Gamma_1:=\left\{\gamma\in C([0,1], H^1(\R^N)): \gamma(0)=0, I_{1,\lambda}(\gamma(1))<0.\right\}.
\eeq
Furthermore, by Corollary \ref{cro:nb-c1},
\beq\lab{eq:20210917-xe8}
\lim_{\lambda\rightarrow 0^+}\frac{m_{1,\lambda}}{\lambda}=0.
\eeq
Similarly, define
\beq\lab{eq:20210917-xe9}
m_{\lambda}=\inf_{\gamma\in \Gamma}\max_{t\in [0,1]}I_{\lambda}(\gamma(t)),
\eeq
where
\beq\lab{eq:20210917-xe10}
\Gamma:=\left\{\gamma\in C([0,1], H^1(\R^N)): \gamma(0)=0, I_{\lambda}(\gamma(1))<0.\right\}.
\eeq
Then, similar to Lemma \ref{lemma:nb-l1}, under the assumption (V) we can prove that
\beq\lab{eq:20210917-xe11}
0<m_\lambda<m_{1,\lambda}.
\eeq
Precisely, $V_0>0$ plays a role in guaranteeing $m_\lambda>0$ and $V_0<V_\infty$ is used to prove that $m_\lambda<m_{1,\lambda}$.
Hence, by \eqref{eq:20210917-xe8} and \eqref{eq:20210917-xe11}, we obtain that
\beq\lab{eq:20210917-xe12}
\lim_{\lambda\rightarrow 0^+}\frac{m_\lambda}{\lambda}=0.
\eeq
Under the assumption (V), for $\lambda\in (0,\lambda_1)$, it is standard to prove the existence of a mountain pass solution $u_\lambda$ to \eqref{eq:20210917-ze1}.
Then, for $\lambda>0$ small,
$$m_\lambda\geq C \lambda \|u_\lambda\|_{H^1(\R^N)}^{2}$$
and \eqref{eq:20210917-xe12} implies that
$$u_\lambda\rightarrow 0~\hbox{in}~H^1(\R^N)~\hbox{as}~\lambda\rightarrow 0^+.$$
Therefore, $\lambda=0$ is also a bifurcation point for \eqref{eq:20210917-ze1}. The proof is complete.

\vskip 0.2in
\noindent
{\bf Acknowledgements.} This work was partially carried out during the second author X.~Zhong visited the Center for Mathematical Sciences (CMS) at Wuhan University of Technology (WUT), he would like to thank CMS-WUT for its hospitality. This work was supported by the Natural Science Foundation of China (No. 11801581, 11871123, 11931012), Guangdong Province Science Foundation (No. 2021A1515010034, 2018A030310082), Guangzhou Science Foundation (No. 202102020225), Chongqing Science Foundation(No. JDDSTD201802) and Chongqing University Science Foundation (No. CXQT21021).


\begin{thebibliography}{100}
\bibitem{ArmstrongSirakov2011}
S.~N.~Armstrong and  B.~Sirakov:
  \newblock Nonexistence of positive supersolutions of elliptic equations via the maximum principle.
  {\it Comm. Partial Differential Equations}, {\bf 36}(11),2011-2047,2011.

\bibitem{BenciFortunato1981}
 V.~Benci and D.~Fortunato:
  \newblock Does bifurcation from the essential spectrum occur?
  {\it Comm. Partial Differential Equations}, {\bf 6}(3),249-272,1981.

\bibitem{BenciFortunato1981-b}
 V.~Benci and D.~Fortunato:
 \newblock Bifurcation from the essential spectrum for odd variational operators.
 {\it Confer. Sem. Mat. Univ. Bari} No.178, 1981.

\bibitem{BerestyckiLions1980}
H.~Berestycki  and P.~L.~Lions:
\newblock Existence of stationary states in nonlinear scalar field equations, Bifurcation phenomena in mathematical physics and related topics, editors C.~Bardos and D.~Br\'ezis, Reidel, 1980.

\bibitem{BerestyckiLions1983}
 H.~Berestycki and P.-L.~Lions:
 \newblock Nonlinear scalar field equations. I. Existence of a ground state.
 {\it Arch. Rational Mech. Anal.} {\bf 82}(4), 313-345, 1983.

\bibitem{Berger1972}
M.~S.~Berger:
\newblock On the existence and structure of stationary states
for a nonlinear Klein-Gordon equation,
{\it J. Functional Anal.}, {\bf 9},249-261,1972.

\bibitem{BHK1983}
A.~Bongers,H.~P.~Heinz and T.~K\"upper:
\newblock Existence and bifurcation theorems for nonlinear elliptic eigenvalue problems on unbounded domains.
{\it J. Differential Equations}, {\bf 47}(3), 327-357,1983.

\bibitem{CLW2017}
S.~W.~Chen, Z.~L.~Liu and Z.~Q.~Wang: A variant of Clark's theorem and its applications for nonsmooth functionals without the Palais-Smale condition. {\it SIAM J. Math. Anal.} {\bf 49}(1), 446--470, 2017.

\bibitem{ChenLi1991}
W.~X.~ Chen and C.~M.~Li:
\newblock{Classification of solutions of some nonlinear elliptic equations}.
{\it Duke Math. J.}, {\bf 63}(3), 615-622,1991.


\bibitem{Clark1973}
D.~ C. ~Clark: A variant of the Lusternik-Schnirelman theory, {\it Indiana Univ. Math. J.}, {\bf 22},65--74,1973.


\bibitem{GidasSpruck1981}
B.~Gidas and J.~Spruck:
\newblock Global and local behavior of positive solutions of nonlinear elliptic equations. {\it Comm. Pure Appl. Math.}, {\bf 35}, 525-598, 1981.

\bibitem{GilbargTrudinger1998}
D.~Gilbarg and N.~S.~Trudinger:
 \newblock Elliptic partial differential equations of second order. Reprint of the 1998 edition.

\bibitem{GuoSeiringer2014}
Y. J. Guo, R. Seiringer: On the mass concentration for Bose-Einstein condensates with attractive interactions, {\it Lett. Math. Phys.}, {\bf 104}, 141-156, 2014.

\bibitem{GuoZengZhou2016}
Y. J. Guo, X. Y. Zeng, H. S. Zhou:  Energy estimates and symmetry breaking in
attractive Bose-Einstein condensates with ring-shaped potentials, {\it Ann. Inst. Henri
Poincar\'e, Anal. Non Lin\'eaire}, {\bf 33}, 809-828, 2016.

\bibitem{HeZou2020}
X.M. He, W.M. Zou:  Bifurcation and multiplicity of positive solutions for nonhomogeneous fractional
Schr\"odinger equations with critical growth. {\it Sci China Math.}, {\bf 63}, 1571-1612, 2020.

\bibitem{Heinz1986}
 H.~P.~Heinz:
  \newblock Nodal properties and bifurcation from the essential spectrum for a class of nonlinear Sturm-Liouville problems.
  {\it J. Differential Equations}, {\bf 64}(1), 79-108,1986.

\bibitem{Heinz1987}
 H.~P.~Heinz:
\newblock Free Ljusternik-Schnirelmantheory and the bifurcation diagrams of certain singular nonlinear systems.
{\it J. Differential Equations}, {\bf 66}, 263--300, 1987.


\bibitem{Jeanjean1999}
L.~Jeanjean:
\newblock Local conditions insuring bifurcation from the continuous spectrum.
\newblock {\em Math. Z.}, {\bf 232}(4), 651-664, 1999.


\bibitem{JeanTanaka2002}
L.~Jeanjean and K.~ Tanaka:
\newblock A Remark on least energy solutions in RN.
{\it Proc. Amer. Math. Soc.}, {\bf 131}, 2399-2408, 2002.

\bibitem{JeanTanaka2003}
 L.~Jeanjean and K.~ Tanaka:
\newblock  A note on a mountain pass characterization of least energy solutions.
{\it Adv. Nonlinear Stud.}, {\bf 3}(4), 445-455, 2003.

\bibitem{JeanZhangZhong2021}
 L.~Jeanjean, J.~J.~Zhang and X.~X.~Zhong:
 \newblock A global branch approach to normalized solutions for the  Schr\"odinger equation.
 \newblock {\em arXiv e-prints}, 2021.\url{https://arxiv.org/abs/2112.05869}

\bibitem{LiYangZou2020}
H.W. Li, Z. Yang, W.M. Zou: Normalized solutions for nonlinear Schr\"odinger equations,
{\it Sci China Math.}, {\bf 50}(8), 1023-1044, 2020. (in Chinese)

\bibitem{LiuWang2015}
Z.~L.~Liu and Z.~Q.~Wang: On Clark's theorem and its applications to partially sublinear problems. {\it Ann. Inst. H. Poincar\'e C Anal. Non Lin\'eaire} {\bf 32}(5), 1015--1037,2015.


\bibitem{QuittnerSouplet2007}
P. Quittner and  P. Souplet:  Superlinear Parabolic Problems. Blow-up, Global Existence and Steady States. Birkh\"auser Advanced Texts: Basler Lehrb\"ucher. Birkh\"auser, Basel (2007).

\bibitem{DautrayLions1984}
 R. Dautray and J.L.Lions: Analyse math\'ematique et calcul num\'erique pour les sciences et les techniques, Masson, Paris, 1984.

\bibitem{Shibata2014}
M.~Shibata:
 \newblock Stable standing waves of nonlinear Schr\"odinger equations with a general nonlinear term. {\it Manuscripta Math.}  {\bf 143}(1-2),221--237,2014.


\bibitem{Strauss1977}
W.~Strauss: Existence of solitaroy waves in higher dimensions,
{\it Comm. Math. Phys.}, {\bf 55}, 149-162, 1977.

\bibitem{Stuart1980}
C.~A.~Stuart:
\newblock  Bifurcation for variational problems when the linearisation has no eigenvalues. {\it J. Functional Analysis}, {\bf 38}(2),169-187, 1980.

\bibitem{Stuart1981}
C.~A.~Stuart:
\newblock  Bifurcation from the continuous spectrum in the $L^2$-theory of elliptic equations on $\R^N$. {\it In Recent Methods in Nonlinear Analysis and Applications}. Liguori, Napoli, 231-300,1981.

\bibitem{Stuart1982}
C.~A.~Stuart:
\newblock  Bifurcation for Dirichlet problems without eigenvalues.
{\it Proc. London Math. Soc.},{\bf 45}(3), 169-192, 1982.

\bibitem{Stuart1983}
C.~A.~Stuart:
\newblock Bifurcation from the essential spectrum.  Equadiff 82 (W\"urzburg, 1982), 575-596, {\it Lecture Notes in Math.}, 1017, Springer, Berlin, 1983.

\bibitem{Stuart1989}
C.~A.~Stuart:
\newblock Bifurcation from the essential spectrum for some noncompact nonlinearities.
{\it Math. Methods Appl. Sci.}  {\bf11}(4), 525-542, 1989.


\bibitem{Wang2001}
Z.~Q.~Wang: Nonlinear boundary value problems with concave nonlinearities near the origin. {\it NoDEA Nonlinear Differential Equations Appl.} {\bf 8}(1), 15--33,2001.



\bibitem{Willem1996}
M.Willem: Minimax theorems. Progress in Nonlinear Differential Equations and their Applications, {\bf 24}. {\it Birkh\"auser Boston, Inc., Boston, MA}, 1996.

\bibitem{YangYang2022}  J.F. Yang, J.G. Yang: Normalized solutions and mass concentration for supercritical nonlinear Schr\"odinger
equations. {\it Sci China Math.}, {\bf 65}, 1383-1412, 2022.

\bibitem{ZhuZhou1988}
X.~P.~Zhu and  H.~S.~Zhou:
\newblock Bifurcation from the essential spectrum of superlinear elliptic equations.
{\it Appl. Anal.},  {\bf28}, 51-66, 1988.

\end{thebibliography}
\end{document}